\documentclass{article}
\usepackage{ijcai19}
\usepackage{times}
\usepackage{soul}
\usepackage{url}
\usepackage[hidelinks]{hyperref}
\usepackage[utf8]{inputenc}
\usepackage[small]{caption}
\usepackage{graphicx}
\usepackage{amsmath}
\usepackage{booktabs}
\usepackage{algorithm}
\usepackage{algorithmic}
\title{Low-Order Model of Biological Neural Networks}

\author{
    Huachuan Wang and James Ting-Ho Lo
    \affiliations
   Department of Mathematics and Statistics, University of Maryland Baltimore County, USA. 
}

\begin{document}

\maketitle

\begin{abstract}
 A biologically plausible low-order model (LOM)
of biological neural networks is a recurrent hierarchical
network of models of dendritic nodes/trees, spiking/nonspiking neurons,
unsupervised/supervised covariance/accumulative learning mechanisms, feedback
connections, and a scheme for maximal generalization. These component models
are motivated and necessitated by making LOM learn and retrieve easily without
differentiation, optimization or iteration, and cluster, detect and recognize
multiple/hierarchical corrupted, distorted and occluded temporal and spatial patterns.
A masking matrix for a dendritic tree, whose upper part comprises model
dendritic encoders, enables maximal generalization on corrupted, distorted and
occluded data. It is a mathematical organization and idealization of dendritic
trees with overlapped and nested input vectors.
A model nonspiking neuron transmits inhibitory graded signals to modulate its
neighboring model spiking neurons. Model spiking neurons evaluate the
subjective probability distribution (SPD) of the labels of the inputs to model
dendritic encoders, and generate spike trains with such SPDs as firing rates.
Feedback connections from the same or higher layers with different numbers of
unit-delay devices reflect different signal traveling times, enabling LOM to
fully utilize temporally and spatially associated information.
Biological plausibility of the component models is discussed. Numerical
examples are given to demonstrate how LOM operates in retrieving,
generalizing, and unsupervised/supervised learning.
\end{abstract}

\section{Introduction}

A learning machine, called a temporal hierarchical probabilistic associative
memory (THPAM), was recently reported \cite{Lo10cody}. The goal to achieve in
the construction of THPAM was to develop a learning machine that learns, with
or without supervision, and retrieves easily without differentiation,
optimization or iteration; and recognizes corrupted, distorted and occluded
temporal and spatial information. In the process to achieve the goal,
mathematical necessity took precedence over biological plausibility. This
top-down approach focused first on minimum mathematical structures and
operations that are required for an effective learning machine with the
mentioned properties.

THPAM turned out to be a functional model of neural networks with many unique
features that well-known models such as the recurrent multilayer perceptron
\cite{HechtN90,PrEuLe00,Bishop06,Haykin09}, associative memories
\cite{Kohonen88,WiBuLo69,Nagano72,Amari89,Suther92,TurAus97}, spiking neural
networks \cite{MaaBis98,GerKis02}, and cortical circuit models
\cite{Martin02,Grange06,Grossb07,GeoHaw09} do not have.

These unique features indicated that THPAM might contain clues for
understanding the operations and structures of biological neural networks. The
components of THPAM were then examined from the biological point of view with
the purpose of constructing a model of biological neural networks with
biologically plausible component models. The components of THPAM were
identified with those of biological neural networks and reconstructed, if
necessary, into biologically plausible models of the same.

This effort resulted in a low-order model (LOM) of biological neural networks.
LOM is a recurrent hierarchical network of biologically plausible models of
dendritic nodes and trees, synapses, spiking and nonspiking neurons,
unsupervised and supervised learning mechanisms, a retrieving mechanism, a
generalization scheme, and feedback connections with delays of different
durations. All of these biologically plausible component models, except the
generalization scheme and the feedback connections, are significantly
different from their corresponding components in THPAM. More will be said
about the differences.

Note that although a dendrite or axon is a part of a neuron, and a
dendro-dendritic synapse is a part of a dendrite (thus a part of a neuron),
they are treated, for simplicity, as if they were separate entities, and the
word \textquotedblleft neuron\textquotedblright\ refers essentially to the
soma of a neuron in this paper.

A basic approximation made in the modeling effort reported here is that LOM is
a discrete-time model and all the spike trains running through it are
Bernoulli processes. The discrete-time approximation is frequently made in
neuroscience. Mathematically, Bernoulli processes are discrete-time
approximation of Poisson processes, which are usually used to model spike
trains. The discrete-time assumption seems similar to that made in the theory
of discrete-time dynamical systems as exemplified by the standard
time-discretization of a differential equation into a difference equation.
However, it is appropriate to point out that the discrete time for LOM and
Bernoulli processes is perhaps more than simply a mathematical approximation.
\begin{enumerate}
\item A spike (or action potential) is not allowed to start within the
refractory period of another, violating a basic property of Poisson processes.
Consecutive spikes in the brain cannot be arbitrarily close and are separated
at least by the duration of a spike including its refractory period, setting
the minimum time length between two consecutive time points in a time line.
\item Neuron or spike synchronization has been discovered in biological
neural networks \cite{Malsbu81,SinGra95}. Such synchronization may originate
from or be driven by synchronous spikes from sensory neurons. 
\item If a
neuron and its synapses integrate multiple spikes in response to sensory
stimuli being held constant (e.g., an image fixated by the retina for about
1/3 of a second), and the neuron generates multiple spikes in the process of
each of such integrations; different time scales and their match-up can
probably be reconciled for the discrete-time LOM and Bernoulli processes.
Before more can be said, the discrete-time approximation is looked upon as
part of the low-order approximation by LOM.
\end{enumerate}
\section{Brief Description of LOM}

LOM is a layered network of processing units (PUs), which is
actually rather simple in structure and operation. Each PU is composed of a
number of\ (model) spiking neurons that jointly output a binary estimate of
the label of the image appearing in the receptive fiels of the PU. Each
neuron contains dendrites, synapses, spiking/nonspiking somas, Hebbian
learning mechanisms, and generalization schemes, which are of course all
computational models.

As a image appears in the receptive field of a PU, the image is encoded by
the dendrites into a dendritic code. If the outer product of this code and
the label estimate output from the PU is added to the memories stored in the
synapses completing the unsupervised Hebbian learning of the image input to
the PU, the PU is called a UPU (unsupervised PU).\ If the label estimated is
replaced with a label provided from outside of the PU in learning the input
image, the supervised Hebbian learning is completed and the PU is called an
SPU (supervised PU). Different dendritic codes are orthogonal, allowing us
to decode the dendritic code of the estimated or provided label for
retrieving the same.

\subsection{Encoding inputs to neurons}

Dendritic trees use more than 60\% of the energy consumed by the brain \cite%
{WongRi89}, occupy more than 99\% of the surface of some neurons \cite%
{FoxBar57}, and are the largest component of neural tissue in volume \cite%
{SirGre87}. Yet, dendritic trees are missing in deep learning machines
(including convolutional neural networks), associative memories \cite%
{Kohone88,HinAnd89,Hassou93} and models of cortical circuits \cite%
{Martin02,Grange06,Grossb07,GeoHaw09}, overlooking a large proportion of the
neuronal circuit.

A key feature of LOM is a novel model of the dendritic encoder. A dendritic
encoder that inputs $v_{\tau } =\left[ 
\begin{array}{ccc}
v_{\tau 1} & v_{\tau 2} & v_{\tau 3}%
\end{array}%
\right] ^{\prime }$ encodes it into the dendritic code $\breve{v}_{\tau }$
by \textbf{the standard parity function }$\phi $\textbf{\ }as follows: 
$\breve{v}_{\tau }= $
$[ 0  \quad v_{\tau 1} \quad v_{\tau 2}  \quad $
$\phi \left( v_{\tau 2},v_{\tau 1}\right) $
$v_{\tau 3}  \quad \phi \left( v_{\tau 3},v_{\tau 1}\right) \quad $
$\phi \left( v_{\tau3},v_{\tau 2}\right) $
$\quad \phi \left( v_{\tau 3},v_{\tau 2},v_{\tau 1}\right)] ^{\prime }$

A graph showing the dendritic encoder is given in Figure. 1. For
example, $\left[ 
\begin{array}{ccc}
1 & 0 & 1
\end{array}
\right] ^{\prime }$ and $\left[ 
\begin{array}{ccc}
0 & 1 & 1
\end{array}
\right] ^{\prime }$\ are encoded into the codes $\left[ 
\begin{array}{cccccccc}
0 & 1 & 0 & 1 & 1 & 0 & 1 & 0
\end{array}
\right] ^{\prime }$\ and $\left[ 
\begin{array}{cccccccc}
0 & 0 & 1 & 1 & 1 & 1 & 0 & 0
\end{array}
\right] ^{\prime }$\ respectively.
\begin{figure}[htbp]
	\centering
		\includegraphics[width=6cm]{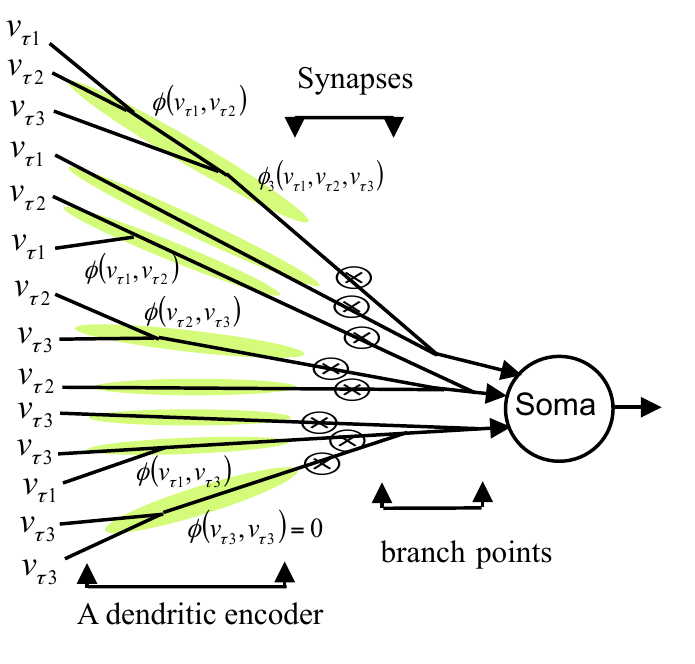}
		\caption{A dendritic encoder.}
\end{figure}
In general, given an input vector $v_{\tau }$ = $\left[ v_{\tau 1},\ldots
,v_{\tau m}\right] ^{\prime }$ to a dendritic encoder, the dendritic code $%
\breve{v}_{\tau }$ of the vector $v_{\tau }$ is
$[ 0 \quad v_{\tau 1}  \quad v_{\tau 2} \quad \phi \left( v_{\tau 2},v_{\tau 1}\right)  \quad $
$v_{\tau 3}  \quad  \phi \left( v_{\tau 3},v_{\tau 1}\right)  \quad  \phi \left( v_{\tau
3},v_{\tau 2}\right) \quad  \phi \left( v_{\tau 3},v_{\tau 2},v_{\tau 1}\right)  \quad  $
$\cdots  \quad \phi \left( v_{\tau m},\cdots ,v_{\tau 1}\right)
]^{\prime }$
Dendritic codes have the orthogonality property: If $v_{\tau }\neq v_{t}$,
then $\left( \breve{v}_{\tau }-\frac{1}{2}\mathbf{1}\right) ^{\prime }\left( 
\breve{v}_{t}-\frac{1}{2}\mathbf{1}\right) =0$. If $v_{\tau }=v_{t}$, then $%
\left( \breve{v}_{\tau }-\frac{1}{2}\mathbf{1}\right) ^{\prime }\left( 
\breve{v}_{t}-\frac{1}{2}\mathbf{1}\right) =2^{\dim v_{\tau }-2}$, where 
\textbf{1 = }$\left[ 
\begin{array}{cccc}
1 & 1 & \cdots & 1%
\end{array}%
\right] ^{\prime }$. This key property is proven in \cite{Lo11neco}. To
avoid curse of dimensionality and enhance generalization capability, only
subvectors of the feature vector $v_{\tau }$\ are expanded as above.

\subsection{Unsupervised/supervised learning and memorizing in synapses}
\begin{figure}[htbp]
	\centering
		\includegraphics[width=8cm]{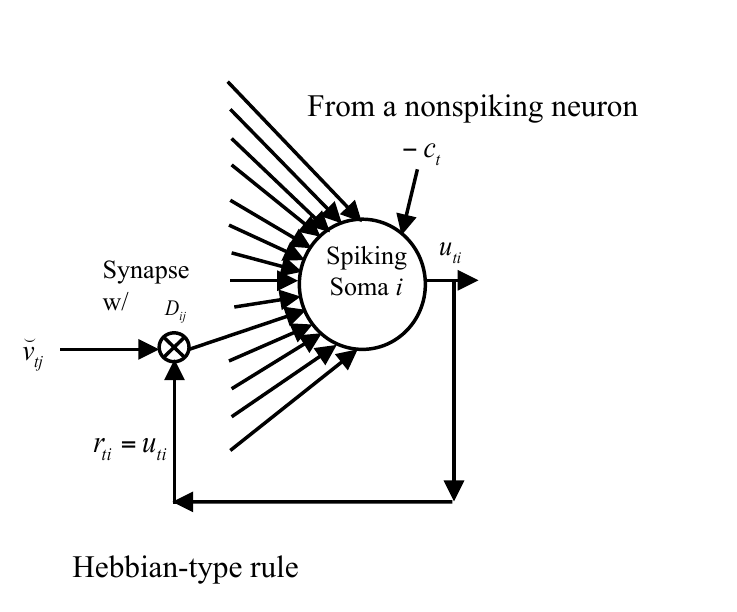}
		\caption{Unsupervised covariance learning.}
\end{figure}
Figure. 2 shows the output $\breve{v}_{tj}$ of the dendritic encoder
going through a synapse with weight $D_{ij}$ represented by $\otimes $ to
reach spiking soma $i$, whose output at time $t$ is $u_{ti}$.The
unsupervised covariance rule that updates the strength $D_{ij}$ of the
synapse receiving $\breve{v}_{tj}$ and feeding spiking neuron $i$ whose
output is $u_{ti}$ follows: 
\begin{equation}
D_{ij}\leftarrow \lambda D_{ij}+\Lambda \left( u_{ti}-\left\langle
u_{ti}\right\rangle \right) \left( \breve{v}_{tj}-\left\langle \breve{v}%
_{tj}\right\rangle \right)  \label{CovUp}
\end{equation}%
where $\Lambda $ is a proportional constant, $\lambda $ is a forgetting
factor that is a positive number less than one, and $\left\langle \breve{v}%
_{tj}\right\rangle $ and $\left\langle u_{ti}\right\rangle $\ denote,
respectively, the average activities of the presynaptic dendritic node $j$
and postsynaptic spiking neuron $i$ over some suitable time intervals.

The outputs $u_{ti}$, $i=1,\ldots ,R,$ of the $R$ spiking somas can be
assembled into a vector, $u_{t}$ = $\left[ 
\begin{array}{cccc}
u_{t1} & u_{t2} & \cdots & u_{tR}%
\end{array}%
\right] ^{\prime }$, and the strengths $D_{ij}$ into a matrix $D$ whose $%
i\times j$-th entry is $D_{ij}$. The vector $u_{t}$ is the label of the
input vector $v_{t}$\ that is selected in accordance with a probability
distribution or membership function by the neurons for unsupervised learning 
\cite{Lo11neco}. Such selection creates a vocabulary for the neurons
themselves.

The synaptic strengths on the connections from the output terminals of a
dendritic encoder to a single nonspiking neuron form a row vector $C$ and
are updated by the unsupervised accumulation rule:%
\begin{equation}
C\leftarrow \lambda C+\frac{\Lambda }{2}\left( \breve{v}_{t}-\left\langle 
\breve{v}_{t}\right\rangle \right) ^{\prime }  \label{CovCup}
\end{equation}

For supervised learning, the output $u_{ti}$ from the spiking soma $i$\ in (%
\ref{CovUp}) is replaced with a component $w_{ti}$ of the label of the input
vector $v_{t}$ provided from outside the LOM. Note that the label is that of
the feature vector in the receptive field of soma $i$.

\subsection{Masking Matrices for Generalization}

If the vector $v_{\tau }$ input to a neuron has not been learned (possibly
due to distortion, corruption and occlusion), it would be ideal if the
largest subvector of $x_{\tau }$ that matches at least one subvector $x_{t}$
stored in the $D$ and $C$ stored in the synapses can be found automatically
and the subjective probability distribution of the label of this largest
subvector of $x_{\tau }$ can be generated. This ideal capability is called
maximal/adjustable generalization. Maximal generalization can actually be
achieved by the use of a masking matrix $M$.

For example, masking the second and fourth components of $\mathbf{1}=\left[ 
\begin{array}{ccccc}
1 & 1 & 1 & 1 & 1%
\end{array}%
\right] ^{\prime }$\ can be done with the multiplication $diag\left[ 
\begin{array}{ccccc}
1 & 0 & 1 & 0 & 1%
\end{array}%
\right] \mathbf{1}$ by the masking matrix $diag\left[ 
\begin{array}{ccccc}
1 & 0 & 1 & 0 & 1%
\end{array}%
\right] $. Masking certain numbers of bits and assigning different weights
to the corresponding masking matrices to deemphasize masking large numbers
of bits allow us to achieve maximal/adjustable generalization. Detailed
description of masking matrices can be found in \cite{Lo11neco}.

\subsection{Estimation of labels by Somas \label{soma}}

Once a vector $v_{\tau }$\ is received and encoded by a dendritic encoder
into $\breve{v}_{\tau }$, $\breve{v}_{\tau }$ is made available to synapses
for learning as well as retrieving of the information about the label of the
input $v_{\tau }$. In response to $\breve{v}_{\tau }$, the masking matrix $M$%
\ computes $M_{jj}\left( \breve{v}_{\tau j}-\left\langle \breve{v}_{\tau
j}\right\rangle \right) $ for all $j$.
\begin{figure}[htbp]
	\centering
		\includegraphics[width=9cm]{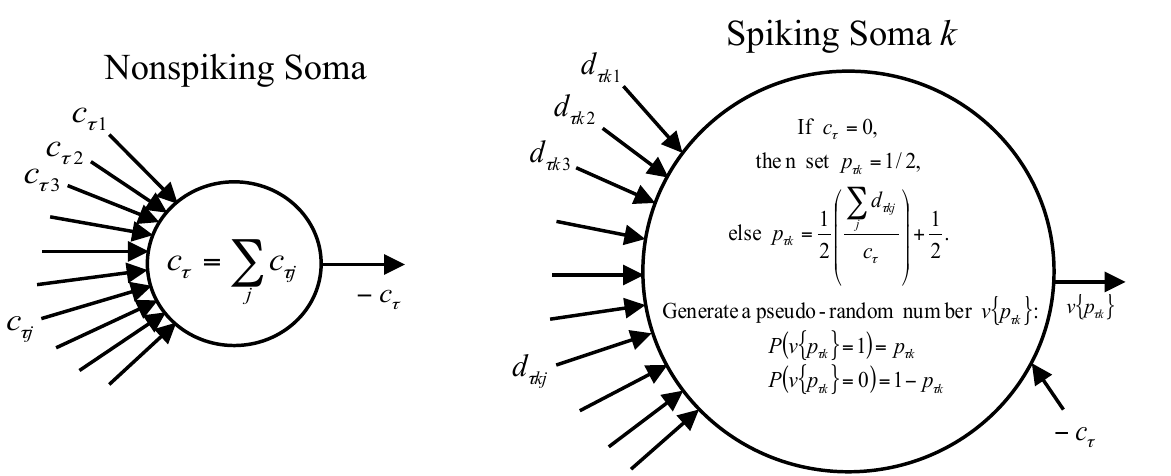}
		\caption{Nonspiking soma and spiking soma k.}
\end{figure}
Synapse $j$ for a nonspiking soma then computes $c_{\tau
j}=C_{j}M_{jj}\left( \breve{v}_{\tau j}-\left\langle \breve{v}_{\tau
j}\right\rangle \right) $, for all $j$, where $M_{jj}$\ is the $j$-th
diagonal entry of $M$.\ As shown in Figure. 3, the model nonspiking soma sums
up $c_{\tau j}$ to obtain the graded signal $c_{\tau }$. Note that the
synaptic weight vector $C$ is defined in (\ref{CovCup}). Because of the
orthogonality property of $\breve{v}_{\tau }$, $t$ = 1, ..., $T$; $c_{\tau }$
is an estimate of the total number of times $\breve{v}_{\tau }$ has been
encoded and stored in $C$. The inhibitory output $-c_{\tau }$is a graded
signal transmitted to each of the mentioned $R$ spiking neurons that
generate a point estimate of the label $r_{t}$ of $v_{\tau }$.

The model spiking soma $k$ is also depicted in Figure. 3. The entries of the $j$%
th row $D_{j}$ of $D$ are the weights or strengths of the synapses for the $%
j $th spiking neuron. In response to $\breve{v}_{\tau }$ produced by the
dendritic encoders, the masking matrix $M$\ and synapses for the $j$th
spiking neuron compute $M_{jj}\left( \breve{v}_{\tau j}-\left\langle \breve{v%
}_{\tau j}\right\rangle \right) $ and $d_{\tau kj}$ = $D_{kj}M_{jj}\left( 
\breve{v}_{\tau j}-\left\langle \breve{v}_{\tau j}\right\rangle \right) $,
respectively. As shown in Figure 3, the $j$th spiking neuron (a model spiking
neuron) sums up $d_{\tau kj}$ to obtain the graded signal $d_{\tau k}$ = $%
\sum_{j}d_{\tau kj}$.

Because of the orthogonality property of $\breve{v}_{\tau }$; $d_{\tau k}$
is an estimate of the total number of times $v_{\tau }$ has been encoded and
stored in $D_{k}$ with the $k$th component $r_{\tau k}$ of $r_{\tau }$ being
1 minus the total number of times $v_{\tau }$ has been encoded and stored in 
$D_{k}$ with the $k$th component $r_{\tau k}$ of $r_{\tau }$ being $0$. The
effects of $M$, $\lambda $ and $\Lambda $ are included in computing said
total numbers, which make $d_{\tau k}$ only an estimate.

Recall that $c_{\tau }$ is an estimate of the total number of times $v_{\tau
}$ has been learned regardless of its labels. Therefore, $\left( c_{\tau
}+d_{\tau k}\right) /2$ is an estimate of the total number of times $v_{\tau
}$ has been encoded and stored with the $k$th component $r_{\tau k}$ of $%
r_{\tau }$ being 1. Consequently, $\left( d_{\tau k}/c_{\tau }+1\right) /2$
is the subjective probability $p_{\tau k}$ that $r_{\tau k}$ is equal to 1.
The $k$th spiking neuron then uses a pseudo-random generator to generate 1
with probability $p_{\tau j}$ and $0$ with probability $1-$ $p_{\tau j}$.
This 1 or $0$ is the output of the $k$th spiking neuron. Biological
justification of the model nonspiking and spiking somas are provided in \cite
{Lo11neco}.

Note that the vector $p_{\tau }=\left[ 
\begin{array}{cccc}
p_{\tau 1} & p_{\tau 2} & \cdots & p_{\tau R}%
\end{array}%
\right] ^{\prime }$ is a representation of a subjective probability
distribution of the label $r_{\tau }$ of the input vector $v_{\tau }$. A
pseudorandom ternary number generator in the $j$th spiking neuron uses $%
p_{\tau j}$ to generate an output denoted by $v\left\{ p_{\tau k}\right\} $
as follows: $v\left\{ p_{\tau k}\right\} $ = 1 with probability $p_{\tau k}$%
, and $v\left\{ p_{\tau k}\right\} $ = $-1$ with probability $1-p_{\tau k}$.
Note also that the outputs of the $R$ spiking neurons in response to $%
v_{\tau }$ form a binary vector $v\left\{ p_{\tau }\right\} $, which is a
point estimate of the label $r_{\tau }$ of $v_{\tau }$.

\subsection{Processing Units (PUs)}
\begin{figure}[htbp]
	\centering
		\includegraphics[width=8.5cm]{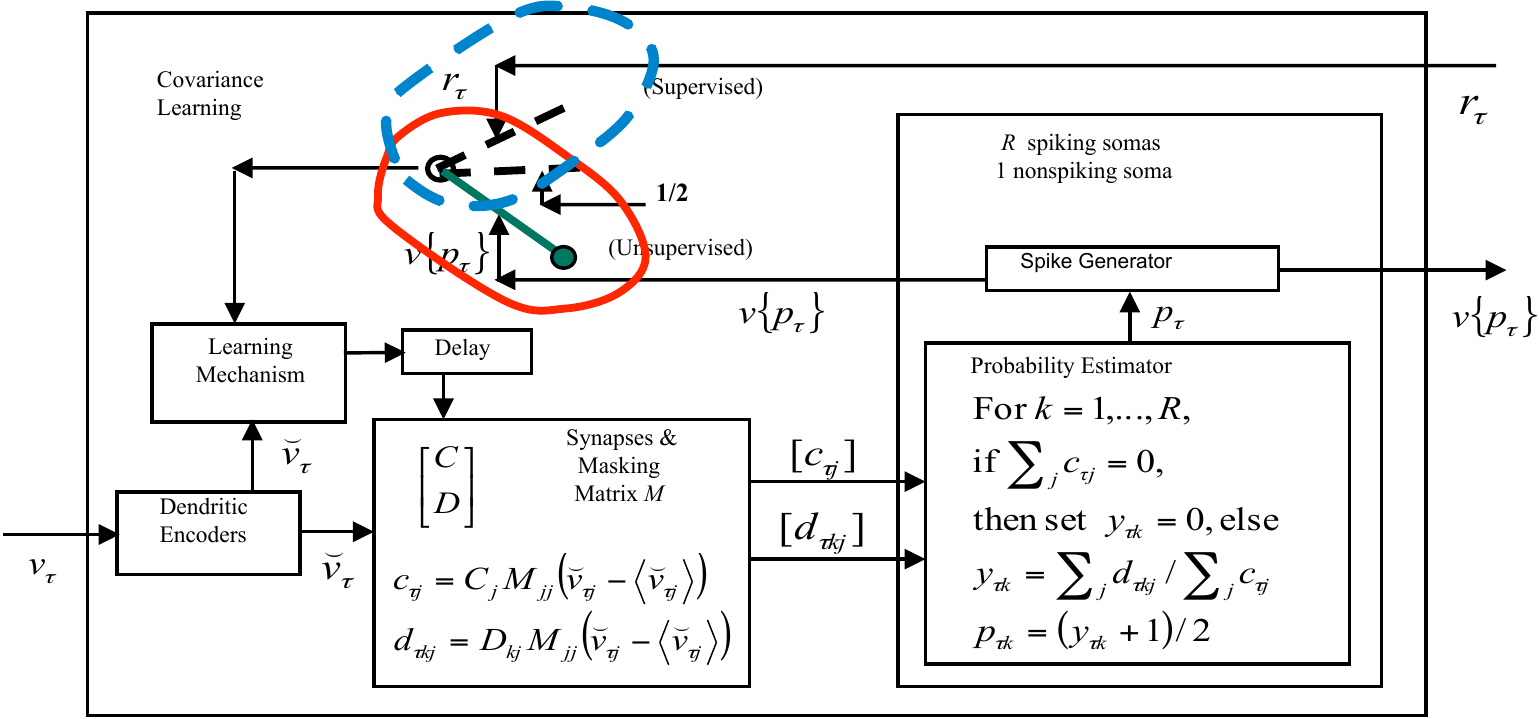}
		\caption{Unsupervised or supervised processing unit (UPU or SPU).}
\end{figure}
LOM organizes a biological neural network into a recurrent network of PUs
(processing units). A schematic diagram of a PU is shown in Figure. 4, which shows how dendritic encoders, synapses, a nonspiking soma, $R$
spiking somas, and learning and retrieving mechanisms are integrated into a
processing unit (PU). The vector $v_{\tau }$ input to a PU is first expanded
by dendritic encoders into a dendritic code $\breve{v}_{\tau }$. $\ \breve{v}%
_{\tau }$ is used to compute $c_{\tau j}$ and $d_{\tau kj}$ using the
synaptic weights $C_{j}$ and $D_{kj}$ respectively. The nonspiking soma in
the PU computes the sum $\sum\limits_{j}c_{\tau j}$, and the $k$th spiking
soma computes $\sum\limits_{j}d_{\tau kj}$ and $p_{\tau k}$ = $\left(
\sum\limits_{j}d_{\tau kj}/\sum\limits_{j}c_{\tau j}+1\right) /2$, which
is the relative frequency that the $k$th digit of the label of $v_{\tau }$
is +1. By a pseudo-random generator, the $k$th spiking soma outputs $%
v\left\{ p_{\tau k}\right\} $, which is +1 with probability $p_{\tau k}$ and
is $-1$ with probability 1 $-$ $p_{\tau k}$, for $k$ = 1, ..., $R$. $%
v\left\{ p_{\tau }\right\} $ = $\left[ 
\begin{array}{ccc}
v\left\{ p_{\tau 1}\right\} & \cdots & v\left\{ p_{\tau R}\right\}%
\end{array}%
\right] ^{\prime }$ is a point estimate of the label of $v_{\tau }$. Note
that use of weights in the masking matrices can facilitate max-pooling of
dendritic encoders in the computation of $p_{\tau k}$ \cite{Lo11neco}.

The green lever circled with the red solid line indicates that the estimated
label $v\left\{ p_{\tau }\right\} $ is used for unsupervised learning. In
this case, the PU is a unsupervised PU (UPU). If the green lever is placed
in the position circled with the blue dashed line, then the handcrafted $%
r_{\tau }$\ is used for supervised learning, and the PU is a supervised PU
(SPU). UPUs in the lowest layer cluster and recognize the lowest level of
pattern elements such as variants of a hyphen, a pipe, a slash, a back
slash, and so on. These pattern elements are integrated from layer to layer
into larger and larger pattern elements and patterns. As long as inputs are
provided to an UPU by sensors or other parts of the LOM, the UPU learns
(without supervision).

By the maximal/adjustable generalization capability (or more specifically,
masking matrix $M$), each UPU acts as a cluster of its input vectors $%
v_{\tau }$. If $v_{\tau }$ or a close version has not been learned by an
UPU, The UPU generates the label of $v_{\tau }$ at random. This enables the
UPU to act as a pattern recognizer by itself. Whenever a handcrafted label $%
r_{\tau }$ is available to an SPU, the SPU learns its input vector $v_{\tau
} $ with $r_{\tau }$. By the maximal/adjustable generalization capability,
the entire cluster(s) constructed by the UPU(s) that provide $v_{\tau }$ is
assigned the same label $r_{\tau }$. This minimizes the amount of
handcrafted labels required.

\subsection{Clustering and Interpreting}

The version of LOM proposed herein consists of a hierarchical network of
UPUs (unsupervised processing units) with feedbacks connections, acting as
pattern recognizers, and a number of offshoot SPUs (supervised processing
units), translating the self-generated labels from UPUs into human language,
which are called the clusterer and interpreter, respectively. An example
clusterer in its entirety for clustering spatial and temporal data is shown
in Figure. 5. The feedback connections in the clusterer have delay devices of
different durations make LOM suitable for recognizing temporal patterns in
video and movie. However, they will not be used in the proposed project.

Once an exogenous feature vector is input to the clusterer, the UPUs perform
retrieving and/or learning from layer to layer starting with layer 1, the
lowest layer. After the UPUs in the highest layer complete performing their
functions, the clusterer is said to have completed one round of retrievings
and/or learnings (or memory adjustments). For each exogenous feature vector,
the clusterer will continue to complete a certain number of rounds of
retrievings and/or learnings.
\begin{figure}[htbp]
	\centering
		\includegraphics[width=8.5cm]{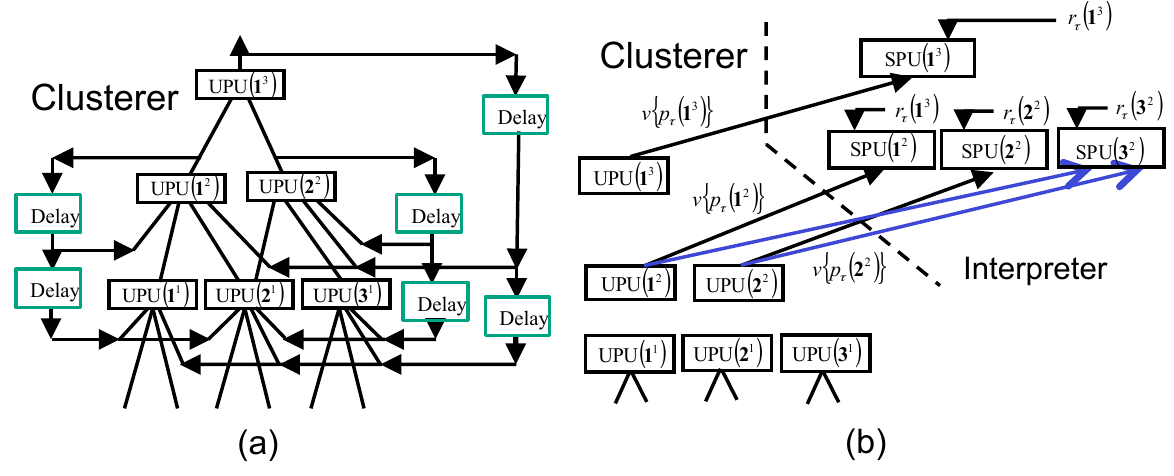}
		\caption{(a)A network of unsupervised processing units(UPUs), (b) offshoot supervised processing units(SPUs)}
\end{figure}
The clusterer in Figure. 5 (a) is also shown in Figure. 5 (b) with the connections and
delay devices removed. The three UPUs in the lowest layer of the clusterer
do not branch out, but each of the three UPUs in the second and third layers
branches out to an SPU. UPU$\left( \mathbf{1}^{2}\right) $ and UPU$\left( 
\mathbf{2}^{2}\right) $ in the second layer have feedforward connections to
SPU$\left( \mathbf{1}^{2}\right) $ and SPU$\left( \mathbf{2}^{2}\right) $
respectively, and UPU$\left( \mathbf{1}^{3}\right) $ in the third layer has
feedforward connections to SPU$\left( \mathbf{1}^{3}\right) $.

The labels, $r_{\tau }\left( \mathbf{1}^{2}\right) $, $r_{\tau }\left( 
\mathbf{2}^{2}\right) $ and $r_{\tau }\left( \mathbf{1}^{3}\right) $, which
are used for supervised learning are provided by the human trainer of the
LOM.

\section{Preliminary Numerical Tests of LOM \label{preliminary}}

\subsection{A simplified supervised learning architecture based on LOM}

This section describes in detail the architecture of a simplified LOM supervised learning neural network model. This architecture comprises four layers, including the input layer, output layer, and two hidden layers. The inputs are $28\times28$ pixel images from the standard MNIST database. In the input layer, a $28\times28$ pixel image is split into 22 per row and 22 per column sub-images, the total number of sub-images is $22\times22 = 484$. The sub-images are generated by sliding an $8\times 8$ pixel window one row or one column each time. An input prepared for a first hidden layer PU is a manually selected16-pixel pattern from the $8\times 8$ pixel sub-image. The selection pattern for each first hidden layer PU is identical. Each PU in each layer is fixed and is always focusing on the same receptive field of the images during the training and testing process.

The input pixels' values are normalized so that the white $(0)$ corresponds to a grayscale value of smaller than 35, and the black $(1)$ corresponds to a grayscale value of larger or equal to 35, which fits the learning scheme. The second hidden layer consists $11\times 11 = 121$ PUs. The input of a PU is formed by combining the output of four lower hidden layer PUs. The lower layer's four PUs are the closest neighbors form a square without overlapping each other.  In this case, the receptive field of a PU in the current layer is extended to $16 \times16 $. 

Within the hidden layers, supervised learning is applied.
In the first hidden layer, before a pattern is prepared as an input
to be learned by the hidden layer PU, the PU will first check
whether the neurons have seen this pattern. If the pattern
has not been learned before, this pattern will be learned using
the image's 4-bit binary label, translated from the
original 0-9 decimal label. However, if the pattern has been
seen before, the neurons can retrieve the pattern's label
by translating probabilities into a label. This pre-check garnered
one pattern that does not have two or more labels, which
greatly confuse the PU when retrieving. This schema
is based on the that the digits share many similar parts. We do
not need to distinguish each part in every digit. What we really
care about is that the pattern is represented following the right
way. For example, if a digit can be covered by the receptive
field of four square PUs in the first hidden layer, because the
top arch pattern of the $"2"$ and $"3"$ are similar, the label retrieved
from the top two PUs of the four PUs will tell the PU
of the next layer this digit is an arch pattern, but they are not
sure whether it is from $"2"$ or $"3"$. After gathering all the four PUs, the PU of the next layer can finally give a prediction based on the patterns shown in each part of the digit.

The second hidden layer does not apply the same learning scheme. In the learning process, a decimal label is translated as a 10-bit one-hot label. Whether the input has been learned or not, it will be learned with the label of the current image. That implies that the same pattern might have two or more labels assigned to it. In the testing process, the 10-bit probability vectors will be treated as the output of the PU.

In the final output layer, each second hidden layer's prediction was analyzed and summarized to produce a final prediction. Each probability vectors will check whether the biggest probability within the vector is larger than 0.85. If the vector's largest probability is smaller than 0.85, that means this pattern might have been assigned many labels and will produce ambiguities towards the final prediction. To obtain a more reliable result, only the vectors with the biggest probability larger than 0.85 will be summed up in the final layer. Finally, the biggest probability after summation will decide the result.

\subsection{Tests on the MNIST dataset}
\begin{figure}[htbp]
	\centering
		\includegraphics[width=8.5cm]{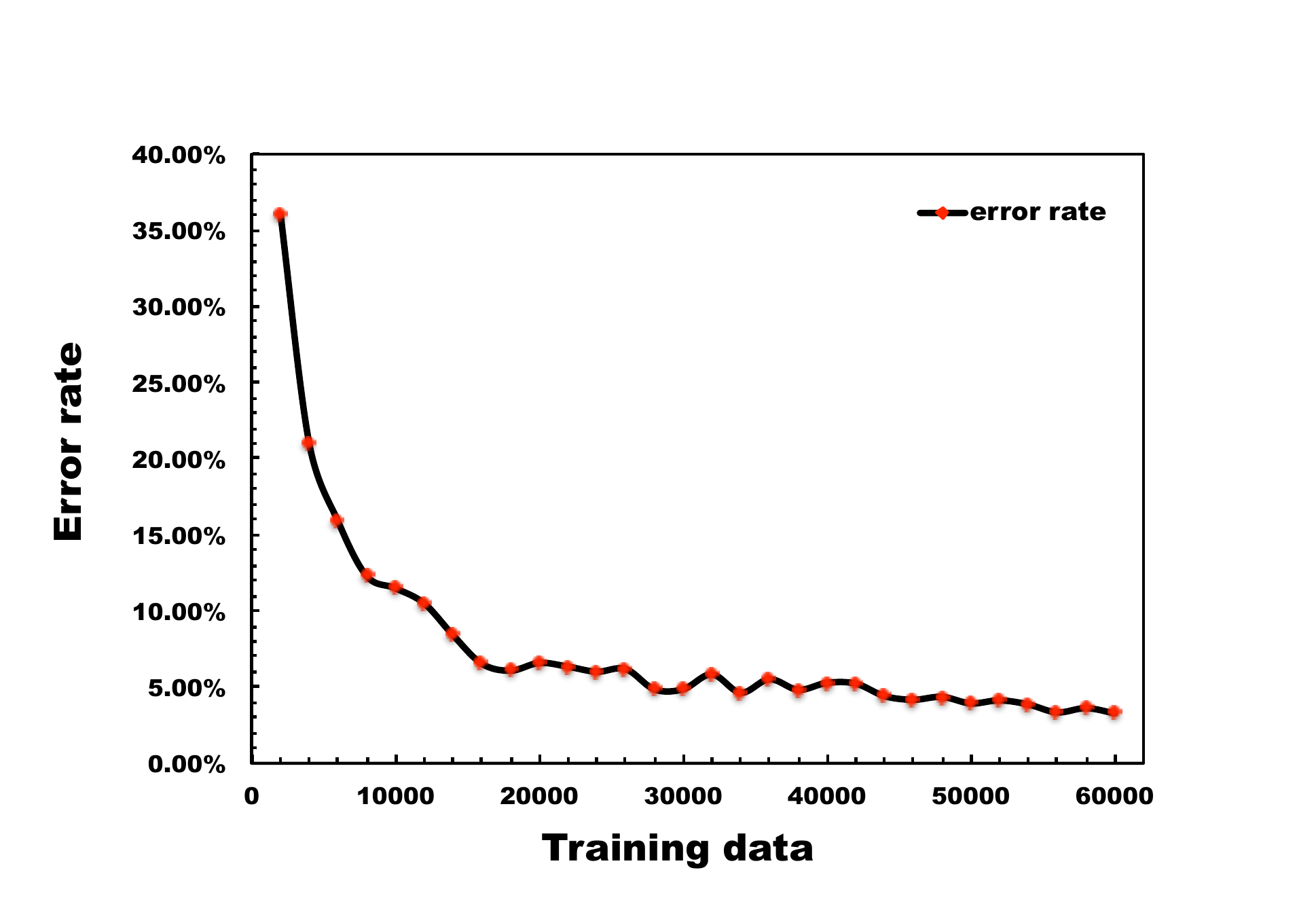}
		\caption{Error rate of real-time learning.}
\end{figure}
The above architecture is developed and tested on the standard.
MNIST dataset consists of 60,000 images of handwritten numerals for training and 10,000 images for testing, each image containing 28$\times $28 pixels. The 60,000 images are divided into 30 bins. Each bin consists of 2000 images. The real-time learning machine-learned each bin sequentially. After learning every 2000 images, the learning machine uses the 10,000 testing images to make a prediction. Following this procedure, a total of 30 predictions of the testing data are made. As shown in Figure. 6. the error rate of each prediction is marked as a red dot. There is an overall downward trend. The error rate drops dramatically from 37\% to lower than 10\%, with only learning the first 12,000 training images. The downward trend slowed down afterward. The final error rate is 3.28\%. The error rate continually dropping around 1\% over the learning procedure of the last 10,000 training images confirmed the LOM model's potential and showed room for improvement. We expect that new LOM architectures will soon be found to achieve an accuracy rate close to the best error rate.
It is appropriate to note that MNIST is a fixed collection of images of
handwritten numerals, each with a single given hand-crafted label, on which training a learning machine requires none of the above desirable learning capabilities. MNIST is, therefore, an ideal kind of dataset for deep learning machine to be trained on without limits on the number of training sessions, the number of
times the weights of deep learning machine can be changed, or the training time in each training session. However, this proposal's main objective is to develop highly accurate LOM architectures with capabilities of real-time, photographic, unsupervised, and hierarchical learning.

\section*{Acknowledgments}
This project is supported by National Science Foundation.
\appendix

\bibliographystyle{named}
\bibliography{ijcai19}

\end{document}